\documentclass[]{amsart}
\usepackage{amsmath,amssymb}
\def\C{\mathbb{C}}

\def\N{\mathbb{N}}
\def\Nquer{\overline{N}}
\def\Z{\mathbb{Z}}

\def\Re{{\rm Re}\;}

\def\s{\mathfrak{s}}
\def\a{\texttt{a}}

\def\ende{ \square}
\def\Ende{~$\ende$\par}

\def\be{\begin{equation}}
\def\ee{\end{equation}}

\def\kzu1#1{\buildrel (#1) \over \longrightarrow}
\newcounter{ABS}

\def\ds{\displaystyle}
\def\ts{\textstyle}

\begin{document}

{\Large\bf Reminiscence Of An Open Problem:

Remarks On Nevanlinna's Four-Value-Theorem}\\

{\large Norbert Steinmetz}\\ {\small Institut f\"ur Mathematik, TU
Dortmund, D-44221 Dortmund,
Germany\\
E-mail: stein@math.tu-dortmund.de\bigskip

\bigskip

{\bf Abstract.} The aim of this paper is to describe the origin,
first solutions, further progress, the state of art, and a new
ansatz in the treatment of a problem dating back to the 1920's,
which still has not found a satisfactory solution and deserves to
be better known.\bigskip

}

\section{\bf Introduction}

 In \cite{polya} G.\ {P\'olya} considered the problem
{\it to determine all pairs $(f,g)$ of distinct entire functions
of finite order, such that $f$ and $g$ assume each of the values
$\a_\nu$ $(1\le \nu \le 3)$ at the same points and with the same
multiplicities,} and solved it as follows.

\medskip{\bf Theorem (P\'olya \cite{polya}).} {\it The functions $f$ and $g$
have a common {Picard} value $\a_2={\ts\frac12}(\a_1+\a_3),$ and
satisfy $(f-\a_2)(g-\a_2)=(\a_2-\a_1)(\a_3-\a_2).$}\medskip

A typical example is $f(z)=e^z,$ $g(z)=e^{-z},$ $\a_1=-1,$
$\a_2=0$, and $\a_3=1$.

In particular, it is not possible that $f$ and $g$ assume each of
{\it four} finite values at the same points and with the same
multiplicities. In modern terminology {P\'olya}'s theorem says
that distinct non-constant entire functions of finite order cannot
{\it share} four finite values by counting multiplicities, and may
share three finite values only in some very particular case.

The background for this theorem was as follows: if two
polynomials, $P$ and $Q$, say, assume integer values at the same
points, then the entire functions $f=e^{2\pi iP}$ and $g=e^{2\pi
iQ}$ have the {Picard} value $0$ and share the value $1$. Thus the
problem arises to determine all entire functions $f$ and $g$ of
finite order having {Picard} value $0$ and sharing the value $1$
by counting multiplicities.

\medskip{\bf Theorem (P\'olya [(\footnote{We follow {Nevanlinna}'s paper
\cite{nevanlinna2};  the reference ``G.\ {P\'olya}, Deutsche
Math.-Ver.\ Bd.\ 32, S.\ 16, 1923'' given there is incorrect.})])}
{\it Under these hypotheses either $f=g$ or else $f=1/g$ holds,
hence either $P-Q$ or else $P+Q$ is a constant $($an integer
multiple of $2\pi i)$.}\medskip

Since the {Picard} value $\infty$ is trivially shared by counting
multiplici\-ties, {P\'olya}'s first theorem may be looked at as a
predecessor of what is nowadays known as Four-Value-Theorem due to
{R. Nevanlinna}, while the proof of his second theorem inspired
{Nevanlinna} to apply his theory of meromorphic functions to
so-called {Borel} identities. To describe these results in more
detail we need some notation; familiarity with the standard
notions and results of {Nevanlinna}'s theory of meromorphic
 functions is assumed, see the standard references {Nevanlinna}
\cite{nevanlinna3} and {Hayman} \cite{hayman}.

Given any pair of distinct meromorphic functions $f$ and $g$
sharing the values $\a_\nu$ ($1\le \nu \le q$), we set
 $$T(r)=\max\,\big(T(r,f), T(r,g)\big),$$
and denote by $S(r)$ the usual remainder term satisfying
 $$S(r)=O(\log(r\;T(r)))\quad(r\to\infty)$$
outside some set $E\subset (0,\infty)$ of finite measure. If $f$
and $g$ have finite [lower] order we have $S(r)=O(\log r)$ and $E$
is empty [$S(r_k)=O(\log r_k)$ on some sequence $r_k\to\infty$].
Furthermore,
 $$\Nquer(r,\a_\nu)=\bar n(0,\a_\nu)+\log r\int_1^r[\bar n(t,\a_\nu)-\bar n(0,\a_\nu)]\,\frac{dt}t$$
denotes the (integrated) Nevanlinna counting function of the
sequence of $\a_\nu$-points of $f$ and $g$, each point being
counted simply despite of multiplicities. Then {Nevanlinna}'s
Theorems may be stated as follows:

\medskip{\bf Five-Value-Theorem (Nevanlinna \cite{nevanlinna2}).}
{\it Let $f$ and $g$ be distinct
meromorphic functions sharing the values{\rm (\footnote{Notably
without hypothesis about the multiplicities.})} $\a_\nu$ $(1\le
\nu \le q)$. Then $q\le 4$, and in case $q=4$ the following is
true:
\begin{itemize}
 \item[$({\rm N}_a)$]$\ds T(r,f)=T(r)+S(r)$ and $\ds T(r,g)=T(r)+S(r).$
 \item[$({\rm N}_b)$]$\ds \sum_{\nu=1}^4 \Nquer(r,\a_\nu)=2\,T(r)+S(r)$.
 \item[$({\rm N}_c)$]$\ds\Nquer\Big(r,\frac
 1{f-g}\Big)=2\,T(r)+S(r)$.{\rm(\footnote{To be modified if
 $\a_4=\infty$: $\ds\Nquer\Big(r,\frac
 1{f-g}\Big)+\Nquer(r,\infty)=2\,T(r)+S(r)$.})}
 \item[$({\rm N}_d)$]$\ds
\Nquer\Big(r,\frac 1{f-b}\Big)=T(r)+S(r)$ and
$\ds\Nquer\Big(r,\frac 1{g-b}\Big)=T(r)+S(r)$ \quad$(b\ne
\a_\nu)$.
\end{itemize}}

The typical example is the same as for P\'olya's Theorem.
\medskip

\medskip{\bf Four-Value-Theorem (Nevanlinna \cite{nevanlinna2}).} {\it Let $f$ and
$g$ be non-constant meromorphic functions sharing values $\a_\nu$
$(1\le \nu \le 4)$, but now specifically {\em by counting
multiplicities}. Then $($relabelling the values, if necessary$)$
$$(\a_1,\a_2,\a_3,\a_4)=(f,g,\a_3,\a_4)=-1$$
holds.$($\footnote{$(a,b,c,d)=\ds\frac{(a-c)(b-d)}{(a-d)(b-c)}$
denotes the cross-ratio of the values $a,b,c,d$.}$)$ In
particular, $\a_1$ and $\a_2$ are {Picard} values of $f$ and $g$,
and $g$ is a {M\"obius} transformation of $f$ that fixes $\a_3$
and $\a_4$ and permutes $\a_1$ and $\a_2$.}\medskip

Both theorems were proved by R.\ {Nevanlinna} in 1926. The novelty
and importance of this paper stems from the new and powerful
methods, nowadays called Nevanlinna Theory, rather than the fact
that meromorphic functions of arbitrary order were considered in
contrast to entire functions of finite order.

\section{\bf Progress and Counterexamples}

In \cite{gundersen1} G.\ Gundersen attributed the question,
whether the Four-Value-Theorem also holds {\it without} the
condition {\it by counting multiplicities}, to L.\ Rubel. This
question or problem, however, was already aware to {Nevanlinna},
who wrote in \cite{nevanlinna2}: {\it Es w\"are nun interessant zu
wissen, ob dieses Ergebnis auch dann besteht, wenn die {\em
Multiplizit\"aten} der betreffenden Stellen nicht
be\-r\"ucksichtigt werden. Einige im ersten Paragraphen gewonnene
Er\-gebnisse {\em [here he refers to conditions $({\rm
N}_a)$-$({\rm N}_d)$]} sprechen vielleicht f\"ur die Vermutung
...}(\footnote{\it It would be interesting to know whether this
result remains true regardless {\em multiplicities}. Some of the
results derived in the first section seem to support the
conjecture ...})

G.\ Gundersen was the first to contribute to that problem. He
proved the

\medskip{\bf (3+1)--Theorem (Gundersen \cite{gundersen1}).} {\it The conclusion of the
Four-Value-Theorem remains true if $f$ and $g$ share four values,
at least {\em three} of them by counting multiplicities.}\medskip

In the same paper, however, Gundersen also provided the first
counterexample to Nevanlinna's conjecture and thus destroyed the
hope for a (0+4)-Theo\-rem. It is easily seen (and this is
Gundersen's example) that {\it the functions
\be\label{ggg}f(z)=\frac{e^z+1}{(e^z-1)^2}{\rm ~and~}
 g(z)=\frac{(e^z+1)^2}{8(e^z-1)}\ee
share the values $1, 0, \infty$ and $-1/8$ in the following
manner: $f$ has only simple zeros and $1$-points, and only double
poles and $(-1/8)$-points, while $g$ has only simple poles and
$(-1/8)$-points, but has only double zeros and $1$-points.}

\medskip Four years later, {Gundersen} was again able to relax the
hypothesis on the number of values shared by counting
multiplicities by proving the

\medskip{\bf (2+2)--Theorem (Gundersen \cite{gundersen2,gundersen3}).}
{\it The conclusion of the Four-Value-Theorem
remains true if $f$ and $g$ share four values, at least {\em two}
of them by counting multiplicities.}\medskip

{Gundersen}'s proof contained a gap, which, however, could be
bridged over (see \cite{gundersen3}) by considering the auxiliary
function which E.\ {Mues} discovered in the early 1980's, but did
not make use of until 1987. This function
$$\Psi=\frac{f'g'(f-g)^2}{\prod\limits_{\nu=1}^4(f-\a_\nu)(g-\a_\nu)}$$
(for $\a_4=\infty$ the factor for $\nu=4$ has to be omitted)
contains the complete information: if $f$ and $g$ share the values
$\a_\nu$ ($1\le \nu \le 4$), then $\Psi$ is an {\it entire}
function, which is {\it small} in the sense that
$$T(r,\Psi)=S(r)$$
holds. It is almost trivial to deduce the 5-Theorem from this
condition of smallness, and it will soon be seen that the whole
progress made afterwards depends almost completely on {Mues}'
function $\Psi$.\medskip

Gundersen's example may also be characterised by additional
properties. This was done by M.\ {Reinders} (a student of Mues) in
several directions, based on the following observations for the
functions (\ref{ggg}):
\begin{itemize}
 \item[(a)]  $f(z_0)=-\frac 12
 \Leftrightarrow g(z_0)=\frac 14$.
\item[(b)] For $\nu$ fixed, one of the functions $f$ and $g$ has
only $\a_\nu$-points of order $2$.
\end{itemize}

\medskip{\bf Theorem (Reinders \cite{reinders1,reinders2}).}
{\it Assume that $f$ and $g$ share mutually distinct values
$\a_\nu$ $(1\le\nu\le 4)$, and that one of the following
conditions holds:
\begin{itemize}
 \item[(a)] There exist
values $a,b\ne \a_\nu$ $(1\le \nu\le 4) $ such that $f(z_0)=a$
implies $g(z_0)=b,$ and the conclusion of the Four-Value-Theorem
does not hold.
 \item[(b$_1$\!)] For every $\nu$ the zeros of $(f-\a_\nu)(g-\a_\nu)$
 have multiplicity $3$.
\item[(b$_2\!)$] For every $\nu$ either $f-\a_\nu$ or else
 $g-\a_\nu$ has only zeros of order $2$.
\end{itemize}
Then up to pre-composition with some non-constant entire function
$h$, and post-composition with some {M\"obius} transformation, $f$
and $g$ coincide with the functions in {Gundersen}'s example
$(\ref{ggg})$.}\medskip

We note that conditions (b$_1\!$) and (b$_2\!$) look quite
different (although they turn out {\it a posteriori} to be equal),
since in case (b$_1\!$), even for fixed $\nu$, each of $f$ and $g$
may have double $\a_\nu$-points, while this is not the case in
(b$_2\!$).\medskip

All attempts to prove or disprove a (1+3)--Theorem failed up to
now (January 2011). Before describing further progress we mention
a counterexample that is quite different from {Gundersen}'s and
has a different origin. It arose from a Comptes Rendus note of H.\
{Cartan} \cite{cartan}, where, in modern terminology, the
following was stated: {\it There do {\em not} exist {\em three}
mutually distinct meromorphic functions sharing four values.} The
proof indicated in \cite{cartan}, however, contained a serious
gap, as {Mues} pointed out to the author in the early 1980's. This
time the gap could not be bridged over, and it took more than
three years to find a way leading to the true statement and, by
the way, to a new counterexample characterised by that theorem.

\medskip{\bf Triple Theorem (Steinmetz \cite{steinmetz}).} {\it Suppose that {\em three}
mutually distinct meromorphic functions $f, g, h$ share four
values, $0, 1,\infty$ and $\a$, say. Then
\begin{itemize}
\item $\a\ne -1$ is a third root of $-1$, and \item $w=f(z), g(z),
h(z)$ are solutions to the algebraic equation
 \be\label{algequ} w^3+3[(\bar
 \a+1)u^2(z)+2u(z)]w^2-3[2u^2(z)+(\a+1)u(z)]w-u^3(z)=0.\ee
\end{itemize}
Here $u=v\circ\gamma$, where $\gamma$ is some non-constant entire
function and $v$ some non-constant solution to the differential
equation
 \be\label{deq}(v')^2=4v(v+1)(v-\a).\ee
Conversely, given $\a\ne 0, -1$, the solutions to equation
$(\ref{deq})$ are elliptic functions, and for $\a=\frac 12(1\pm
i\sqrt{3})$ and any non-constant entire function $\gamma$, the
solutions to equation $(\ref{algequ})$ with $u=v\circ\gamma$
provides three meromorphic functions sharing the values $0, 1,
\infty,$ and $\a$.}\medskip

In the most simple case $\gamma(z)=z$ these functions are elliptic
functions of elliptic order six. They share the values $0, 1,
\infty,$ and $\a$ in the following manner: every period
parallelogram contains three $c$-points ($c\in\{0, 1, \infty,
\a\}$), each being simple for two of these functions, and having
multiplicity {\it four} for the third one. Thus the sequence of
$c$-points is divided in a natural way into three subsequences
having asymptotically equal counting functions. \medskip

Reinders constructed a counterexample that is quite different
from {Gundersen}'s. It is well-known that the non-constant
solutions of the differential equation
 \be\label{deqreinders}(u')^2=12\,u(u+1)(u+4),\ee
are elliptic functions of elliptic order two (actually
$u(z)=\wp(\sqrt3 z+c)-5/3$, where $\wp$ is the specific P-function
of Weierstrass that satisfies the differential equation
$(\wp')^2=4(\wp-\frac 53)(\wp-\frac 23)(\wp+\frac 73)$). {\it Then
for any such function,
\be\label{rex}f=\frac1{8\sqrt{3}}\,\frac{uu'}{u+1} {\rm ~and~}
 g=\frac1{8\sqrt{3}}\,\frac{(u+4)u'}{(u+1)^2}\ee
share the values $-1,0,1,\infty$ in the following manner: each of
these values is assumed in an alternating way with multiplicity
$1$ by one of the functions $f$ and $g$, and with multiplicity $3$
by the other one} (e.g., $f$ has triply zeros at the zeros of $u$,
and simple zeros when $u+4=0$.)\medskip

Again this example can be characterised by this particular
property.

\medskip{\bf Theorem (Reinders \cite{reinders3}).}
{\it Let $f$ and $g$ be meromorphic functions sharing finite
values $\a_\nu$ $(1\le\nu\le 4)$, such that $(f-\a_\nu)(g-\a_\nu)$
has always zeros of multiplicity at least $4$. Then up to
pre-compo\-sition with some non-constant entire function and
post-composition with some {M\"obius} transformation, the
functions $f$ and $g$ coincide with those in example
$(\ref{rex})$.}\medskip

{\it Remark.} Suppose that there exist positive integers $p<q$,
such that every shared value $\a_\nu$ is assumed by $f$ with
multiplicity $p$ and by $g$ with multiplicity $q$, or vice versa.
Then considering Mues' function yields $p=1$, and from
$$(1+q)\Nquer(r,\a_\nu)=N\Big(r,\frac1{f-\a_\nu}\Big)+
 N\Big(r,\frac1{g-\a_\nu}\Big)+S(r)\le 2T(r)+S(r)$$
and condition $(N_b)$ easily follows $2(1+q)T(r)\le 8T(r)+S(r),$
hence  $2\le q\le 3$. On combination with the theorems of Reinders
this immediately yields:

\medskip{\bf Theorem (Chen, Chen, Ou \&
Tsai \cite{CCCT}).} {\it Suppose $f$ and $g$ share four values
$\a_\nu$, and assume also that there are integers $1\le p<q$ such
that each $\a_\nu$-point is either a $(p,q)$-fold point for
$(f,g)$ or for $(g,f)$. Then the conclusion of Reinders' theorems
\cite{reinders1,reinders2,reinders3}
hold.}\medskip

\section{\bf Proof of the (2+2)--Theorem}

Any progress till now relies on the fact that ``(2+2) implies
(4+0)'', and was based on {Mues}' auxiliary function technique
\cite{mues}. To show the power of this method we will next give
independent and short proofs of the Four-Value-, the (3+1)-, and
the (2+2)--Theorem.

\medskip{\bf Proof of the (4+0)--Theorem---Mues \cite{mues}}. Suppose $f$ and
$g$ share four values $\a_\nu$ $(1\le\nu\le 4)$ by counting
multiplicities. Then at least two of the values $\a_\nu$ satisfy
$\Nquer(r,\a_\nu)\ne S(r)$. We may assume that $\a_4=\infty$ does,
hence $N(r,f)=\Nquer(r,\infty)+S(r)$ and
$N(r,g)=\Nquer(r,\infty)+S(r)$ hold, while
$N(r,1/f')+N(r,1/g')=S(r)$ is an easy consequence of the strong
assumption ``by counting multiplicities''. Thus the auxiliary
function
 $$\phi=\frac{f''}{f'}-\frac{g''}{g'}$$
satisfies $T(r,\phi)=S(r)$, but vanishes at all poles that are
simple for $f$ and $g$. This yields a contradiction, namely
$\Nquer(r,\infty)=S(r)$, if $\phi\not\equiv 0$. If, however,
$\phi$ vanishes identically, $f=g$ follows at once. \Ende

\medskip{\bf Proof of the (3+1)--Theorem---Rudolph \cite{rudolph}}. Suppose $f$ and $g$
share the finite values $\a_\nu$ $(1\le\nu\le 3)$ by counting
multiplicities, and $\a_4=\infty$ without further hypothesis. Then
 $$\phi=\frac{f'}{\prod\limits_{\nu=1}^3(f-\a_\nu)}-\frac{g'}{\prod\limits_{\nu=1}^3(g-\a_\nu)}$$
satisfies $T(r,\phi)=S(r)$ and vanishes at the poles of $f$ and
$g$. Thus we have either $\phi\not\equiv 0$ and
 $\Nquer(r,\infty)=S(r),$
or else $\phi$ vanishes identically, the latter meaning that also
$\a_4=\infty$ is shared by counting multiplicities. On the other
hand it is not hard to show that the hypotheses
$``\Nquer(r,\infty)=S(r)$'' and {$``\a_4=\infty$ is shared by
counting multiplicities''} are equally strong in the sense that
(on combination with the hypotheses about the values $\a_\nu$)
they lead to the same conclusions. \Ende

\medskip{\bf Proof of the (2+2)--Theorem.} We may assume that the values $0$ and
$\infty$ are shared by counting multiplicities, while the other
values are $\a_1$ and $\a_2=1/\a_1$ (if the latter does not hold
{\it a priori}, we consider $cf$ and $cg$ instead of $f$ and $g$,
with $c$ satisfying $c^2\a_1\a_2=1$). The auxiliary function
$$\phi=
\frac{f''}{f'}+2\frac{f'}{f}-\sum_{\nu=1}^2\frac{f'}{f-\a_\nu}-
\frac{g''}{g'}-2\frac{g'}{g}+\sum_{\nu=1}^2\frac{g'}{g-\a_\nu}$$
is regular at all $\a_\nu$-points ($\nu=1,2$) of $f$ and $g$
despite their multiplicities, and is also regular at the zeros and
poles that are simple for $f$ and $g$. Thus $\phi$ satisfies
$T(r,\phi)=S(r)$, and
 $$\phi(z_\infty)^2=(\a_1+\a_2)^2\Psi(z_\infty)$$
holds at any pole which is simple for $f$ and $g$; here
$\Psi$ 
again is Mues' function. Hence either $\Nquer(r,\infty)=S(r)$ or
else $\phi^2=(\a_1+\a_2)^2\Psi$ holds. Repeating this argument
with $F=1/f$, $G=1/g$ and corresponding function
$$\begin{array}{rcl}
\chi&=&\ds\frac{F''}{F'}+2\frac{F'}{F}-\sum_{\nu=1}^2\frac{F'}{F-1/\a_\nu}-
\frac{G''}{G'}-2\frac{G'}{G}+\sum_{\nu=1}^2\frac{G'}{G-1/\a_\nu}\cr
&=&\ds\frac{f''}{f'}-2\frac{f'}{f}-\sum_{\nu=1}^2\frac{f'}{f-\a_\nu}
 -\frac{g''}{g'}+2\frac{g'}{g}+\sum_{\nu=1}^2\frac{g'}{g-\a_\nu}\end{array}$$
instead of $f$, $g$ and $\phi$, and noting that
$1/\a_1+1/\a_2=\a_1+\a_2$ and $\Psi$ is also Mues' function for
$F$ and $G$, it follows that essentially four possibilities remain
to be discussed:

\begin{itemize}\item [(a)] $\Nquer(r,0)+\Nquer(r,\infty)=S(r)$;
\item [(b)] $\phi=\chi$; \item [(c)] $\phi=-\chi$;\item [(d)]
$\Nquer(r,\infty)=S(r)$ 
and $\chi^2=(\a_1+\a_2)^2\Psi$.
\end{itemize}

{ad (a)---}From (N$_b$) follows
$\Nquer(r,\a_1)+\Nquer(r,\a_2)=2T(r)+S(r)$, hence the sequence of
$\a_\nu$-points ($\nu=1, 2$) which have different multiplicities
for $f$ and $g$ have counting function $S(r)$. The conclusion of
the Four-Value-Theorem follows immediately from Mues' proof of
that theorem.

\medskip{ad (b)---}From
 $\ds\phi-\chi=4\Big(\frac{f'}{f}-\frac{g'}{g}\Big)= 0$
follows $g=cf$, hence $g=f$ if there exists some common
$\a_\nu$-point. Otherwise $\a_1$ and $\a_2$ are Picard values for
$f$ and $g$, hence the Four-Value-Theorem holds.

\medskip{ad (c)---}Here we have
$\ds\frac{f''}{f'}-\sum_{\nu=1}^2\frac{f'}{f-\a_\nu}-\frac{g''}{g'}+
 \sum_{\nu=1}^2\frac{g'}{g-\a_\nu}= 0,$
hence
$$\ds\frac{f'}{(f-\a_1)(f-\a_2)}=
 \kappa\frac{g'}{(g-\a_1)(g-\a_2)}\quad(\kappa\ne 0)$$
and
 $$\ds\frac{f-\a_1}{f-\a_2}=C\Big(\frac{g-\a_1}{g-\a_2}\Big)^\kappa$$
for some $\kappa\in\Z\setminus\{0\}$ and $C\ne 0.$ Then
$\kappa=\pm 1$ follows from $T(r,f)=T(r,g)+S(r)$, and two cases
remain to be discussed: $\kappa=$

\begin{itemize}\item[$-1.$] The values $\a_\nu$ are Picard values for $f$
and $g$, hence $f$ and $g$ share all values by counting
multiplicities. \item[$1.$] From
$\ds\frac{f-\a_1}{f-\a_2}=C\frac{g-\a_1}{g-\a_2}$ follows that $f$
and $g$ also share the values $\a_1$ and $\a_2$, hence all values,
by counting multiplicities.
\end{itemize}

{ad (d)---}Following Mues \cite{mues} we consider the auxiliary
functions
$$\eta_1=\chi-(\a_1+\a_2)\frac{f'(f-g)}{f(g-\a_1)(f-\a_2)}
{\rm~and~}
 \theta_1=\chi-(\a_1+\a_2)\frac{g'(f-g)}{g(f-\a_1)(g-\a_2)}$$
and the corresponding functions $\eta_2$ and $\theta_2$ (with
$\a_1$ and $\a_2$ permuted). It is easily seen that
 $$T(r,\eta_\nu)\le T(r)-\Nquer(r,\a_\nu)+S(r)  {\rm~and~}
 T(r,\theta_\nu)\le T(r)-\Nquer(r,\a_\nu)+S(r)$$
holds (see also \cite{mues}). Now each of these functions
vanishes at the zeros of $f$ and $g$. If the functions $\eta_1$
and $\theta_1$ as well as $\eta_2$ and $\theta_2$ do not
simultaneously vanish identically, we obtain (note that
$\Nquer(r,\infty)=S(r)$)
$$2\Nquer(r,0)\le
 2T(r)-\Nquer(r,\a_1)-\Nquer(r,\a_2)+S(r)=\Nquer(r,0)+S(r),$$
hence $\Nquer(r,0)=S(r)$. So we are back in case (a), and the
(2+2)--Theorem is proved completely in that case.

On the other hand, $\eta_1=\theta_1\equiv 0$, say, and
$\a_1+\a_2\ne 0$ yield
 $$\frac{(f-\a_1)f'}{f(f-\a_2)}=\frac{(g-\a_1)g'}{g(g-\a_2)},$$
hence $f$ and $g$ share the value $\a_1$ by counting
multiplicities, so that the (3+1)--Theorem gives the desired
result.

The case $\a_1+\a_2=0$ (hence $\a_1=i$ and $\a_2=-i$, on combination with $\a_1\a_2=1$) has to be
treated separately. From $\chi\equiv 0$ follows
\be\label{varphi}\frac{f'}{f^2(f-i)(f+i)}=\kappa\frac{g'}{g^2(g-i)(g+i)}\quad(\kappa\ne
 0).\ee
Since the values $\pm i$ are not Picard values for $f$ and $g$
(otherwise we were already in the (3+1)-case), $\kappa$ (or
$1/\kappa$) is a positive integer, hence $f$ assumes the values
$\pm i$ ``always'' with multiplicity $\kappa$, while $g$ has
``only'' simple $\pm i$-points (up to a sequence of points with
counting function $S(r)$). Now (\ref{varphi}) is equivalent to
 $$\frac{f'}{f^2}-\kappa\frac{g'}{g^2}=\frac i2
\Big(\frac{f'}{f-i}-\frac{f'}{f+i}-\kappa\frac{g'}{g-i}+\kappa\frac{g'}{g+i}\Big),$$
and the right hand side, denoted $\varphi$, is regular at $\pm
i$-points and vanishes at simple zeros $z_0$ of $f$ and $g$ (note
that $g'(z_0)=\kappa f'(z_0)$). Also from $T(r,\varphi)=S(r)$
follows either $\Nquer(r,0)=S(r)$, which leads back to case (a),
or else $\varphi\equiv 0$, equivalently
$$\frac{f'}{f^2}-\kappa\frac{g'}{g^2}=0 {\rm ~and~} \frac 1f=\frac
 \kappa g+c$$
with $\pm 1/i=\pm\kappa/i+c$. This is only possible if $\kappa=1$
and $c=0$. \Ende

\section{\bf Progress after {Gundersen}}

As already mentioned, all attempts to prove or disprove a
(1+3)--Theorem failed up to now (January 2011), and so many
authors looked for additional conditions or switched to related
problems---but the latter is not the subject of the present paper.
In \cite{mues} E.\ Mues introduced the quantity
$$\tau(\a_\nu)=\liminf_{r\to\infty~(r\notin E)}
\frac{N_\s(r,\a_\nu)}{\Nquer(r,\a_\nu)}\quad\hbox{\rm
if~~}\Nquer(r,\a_\nu)\ne S(r),$$ and $\tau(\a_\nu)=1$ else; here
$N_\s(r,\a_\nu)$ denotes the counting function of those
$\a_\nu$-points which are {\it $\s$imultaneously $\s$imple} for
$f$ and $g$, and $E$ is the exceptional set for $S(r)$; we note
that $\tau(\a_\nu)=1$ in particular holds if $\a_\nu$ is shared by
counting multiplicities, and also if  $\a_\nu$ is a {Picard} value
for $f$ and $g$.

\medskip{\it Example.} We have $\tau(\a_\nu)=0$ in the counterexamples of
Gundersen and Reinders, and $\tau(\a_\nu)=\frac 13$ for any pair
of the author's triple.\medskip

In her diploma thesis, E.\ {Rudolph} \cite{rudolph} proved some
results in terms of the quantities $\tau(\a_\nu)$ or their natural
generalisations $\ds\tau(\a_{\nu},\a_{\mu})$,
$\ds\tau(\a_{\nu},\a_{\mu},\a_{\kappa})$ and
$\ds\tau(\a_{1},\a_{2},\a_{3},\a_{4})$. Her proofs were based on
the methods developed in the paper \cite{mues}, which appeared
later than \cite{rudolph}, but was written earlier. In the sequel
we shall derive several results, which become much more apparent
when stated as inequalities involving the counting functions
$\Nquer(r,\a_\nu)$ and $N_\s(r,\a_\nu).$ Nevertheless they may be
credited to Rudolph and Mues for the underlying idea.

\medskip{\bf Key Lemma.} {\it Suppose that $f$ and $g$ share
distinct values $\a_\nu$ $(1\le\nu\le 4)$. Then either the
conclusion of the Four-Value-Theorem holds or else the following
is true:
\begin{itemize}
\item[(${\rm R}_a$)]\quad $\frac 32N_\s(r,\a_\kappa)
+N_\s(r,\a_\nu) \le\Nquer(r,\a_\kappa)+\Nquer(r,\a_\nu)+S(r)$ for
$\kappa\ne\nu$;\item[] \quad for $(\a_1,\a_2,\a_3,\a_4)=-1$ the
factor $\frac 32$ can be replaced by $2$. \item[(${\rm
R}_b$)]\quad $\ds N_\s(r,\a_\kappa)+ \sum_{\mu\ne
\lambda}N_\s(r,\a_\mu) \le
\sum_{\mu\ne\lambda}\Nquer(r,\a_\mu)+S(r)$ for $\kappa\ne\lambda.$
\item[(${\rm R}_c$)]\quad $\ds\Nquer(r,\a_\kappa)+\sum_{\mu\ne
\kappa}N_\s(r,\a_\mu)\le \sum_{\mu\ne
\kappa}\Nquer(r,\a_\mu)+S(r)$. \item[(${\rm R}_d)$]\quad
$\ds\sum_{\mu=1}^4
N_\s(r,\a_\mu)+\sum_{\mu\ne\kappa,\nu}N_\s(r,\a_\mu) \le
2\sum_{\mu\ne\kappa,\nu}\Nquer(r,\a_\mu)$ for $\kappa\ne\nu.$
\item[(${\rm R}_e$)]\quad $\ds 4\sum_{\kappa\ne\lambda}
N_\s(r,\a_\kappa) \le 3\sum_{\kappa\ne\lambda}
\Nquer(r,\a_\kappa)+S(r).$ \item[(${\rm R}_f$)]\quad $\ds
3\sum_{\mu=1}^4 N_\s(r,\a_\mu)\le
 2\sum_{\mu=1}^4 \Nquer(r,\a_\mu)+S(r)=4T(r)+S(r)$.\end{itemize}}\medskip

Several more or less recent results obtained by different authors
can easily be derived from the previous lemma. Since, however, the
thesis \cite{rudolph} has never been published, these results can
be looked upon as independent discoveries. The proof of the Key
Lemma will be given in the next section.

\medskip{\bf Corollary.} {\it Suppose that $f$ and $g$ share the values $\a_\nu$
$(1\le\nu\le 4)$. Then the conclusion of the Four-Value-Theorem is
true, provided one of the following hypotheses is assumed in
addition:}
\begin{itemize}
 \item[(A)] {\it {\rm One} value is shared by
counting multiplicities, and some other satisfies
$\tau(\a_\kappa)>\frac 23$; for $(\a_1,\a_2,\a_3,\a_4)=-1$ the
constant $\frac 23$ can be replaced by $\frac 12$}.{\bf---(E.\
Mues \cite{mues}).(\footnote{It is remarkable that Mues did not
refer to Gundersen's (2+2)--Theorem, hence, in particular, gave an
independent proof for it. The other authors actually proved that
the hypotheses of the (2+2)--Theorem follow from their own, thus
their theorems generalising Gundersen's (2+2)-result actually
depend on it.})}
 \item[(B)] {\it {\rm One} value is shared
by counting multiplicities, while the other values satisfy
$\tau(\a_\nu)>\frac 12$}.{\bf ---(J.P.~{Wang} \cite{wang2},
B.~{Huang}~\cite{huang}).}
 \item[(C)]
{\it {\rm One} value is shared by counting multiplicities and
simultaneously satisfies $\Nquer(r,\a_\kappa)\ge ({\ts\frac
45}+\delta)\,{T(r)}$ for some $\delta>0$ on some set of infinite
measure.{\bf ---(G.\ {Gundersen} \cite{gundersen}).}}
 \item[(D)]  {\it {\rm Two} of the values satisfy
$\tau(\a_\nu)>\frac 45$; for $(\a_1,\a_2,\a_3,\a_4)=-1$ the
constant $\frac 45$ can be replaced by $\frac 23$.}{\bf
---(E.\ Rudolph \cite{rudolph}, S.P.\ {Wang}~\cite{wang},
B.~{Huang}~\cite{huang}).} \item[(E)] {\it {\rm Three} of the
values satisfy $\tau(\a_\nu)>\frac 34$.}{\bf
---(E.\ Rudolph \cite{rudolph}, G.\ Song~\& J.~Chang
\cite{song-chang}).} \item[(F)] {\it {\rm All} values satisfy
$\tau(\a_\nu)>\frac 23$}.{\bf
---(E.\ Rudolph \cite{rudolph}, H.\ Ueda \cite{ueda}(\footnote{Ueda
has a slightly stronger result involving some additional term,
which, however, cannot be controlled.}), J.P.~Wang~\cite{wang2},}
also conjectured by {\bf G.\ Song \& J.~Chang \cite{song-chang})}.
\end{itemize}\medskip

{\it Proof.}  Assuming that the conclusion of Nevanlinna's
Four-Value-Theorem does {\it not} hold, we will derive a
contradiction to the respective hypothesis by applying one of the
inequalities (R$_a$) -- (R$_f$).

(A) is true since
 $$\Nquer(r,\a_\nu)=N_\s(r,\a_\nu)+S(r)$$
for some $\nu$ and $({\rm R}_a)$ for the same $\nu$ imply
$${\ts\frac 32}N_\s(r,\a_\kappa)
 \le\Nquer(r,\a_\kappa)+S(r)\quad(\kappa\ne \nu),$$
with $\frac 32$ replaced by $2$ if $(\a_1,\a_2,\a_3,\a_4)=-1$.

To prove (B) we may assume that  $\a_4$ is shared by counting
multipli\-cities, thus
 $$\Nquer(r,\a_4)=N_\s(r,\a_4)+S(r)$$
holds. From $({\rm R}_c)$ with $\kappa<4$ then follows
 $$\Nquer(r,\a_\kappa)+\sum_{\mu\ne\kappa,4}N_\s(r,\a_\mu)\le
 \sum_{\mu\ne\kappa,4}\Nquer(r,\a_\mu)+S(r),$$
 and adding up for $\kappa=1, 2, 3$ yields
 $$\sum_{\kappa=1}^3\Nquer(r,\a_\kappa)+2\sum_{\mu=1}^3
 N_\s(r,\a_\mu)\le2\sum_{\mu=1}^3\Nquer(r,\a_\mu)+S(r).$$

(C) is obtained from (${\rm R}_a$) as follows. We assume that
$\a_\kappa$ is shared by counting multiplicities, hence $\Nquer(r,\a_\kappa)=N_\s(r,\a_\kappa)+S(r)$ holds.
Adding up for $\nu\ne\kappa$ we obtain
 $${\ts\frac 32}\Nquer(r,\a_\kappa)+\sum_{\nu\ne\kappa}N_\s(r,\a_\nu)\le
\sum_{\nu\ne\kappa}\Nquer(r,\a_\nu)+S(r)=2T(r)-\overline
 N(r,\a_\kappa)+S(r).$$

(D) again follows from $({\rm R}_a)$ by adding up the symmetric
inequalities
$$\begin{array}{rcl}
{\ts\frac 32}N_\s(r,\a_\nu)+N_\s(r,\a_\kappa)
&\le&\Nquer(r,\a_\kappa)+\Nquer(r,\a_\nu)+S(r)\cr {\ts\frac
32}N_\s(r,\a_\kappa)+N_\s(r,\a_\nu)
&\le&\Nquer(r,\a_\kappa)+\Nquer(r,\a_\nu)+S(r).\end{array}$$ Again
we note that $\frac 45$ can be replaced by $\frac 23$ if
$(\a_1,\a_2,\a_3,\a_4)=-1$.

Finally, (E) and (F) follow from (${\rm R}_e$) and (${\rm R}_f$), respectively. \Ende

\section{\bf Proof of the Key Lemma}

The proof of the Key Lemma is based on {Mues}' auxiliary function
technique \cite{mues}. As long as no particular hypotheses are
imposed, the proofs and results are symmetric and M\"obius
invariant, this meaning that everything proved for
$\a_1,\a_2,\a_3,\a_4$  holds for arbitrary permutations, and that
three of the four values can be given prescribed numerical values.

\medskip To prove (${\rm R}_a$) we proceed as in the proof of the
(2+2)--Theorem. We set $\kappa=4$, $\a_4=\infty$, and $\nu=3$, and
consider Mues' function $\Psi$ and the auxiliary function
$$\phi=
\frac{f''}{f'}+2\frac{f'}{f-\a_3}-\sum_{\nu=1}^2
\frac{f'}{f-\a_\nu}-\frac{g''}{g'}-2\frac{g'}{g-\a_3}
 +\sum_{\nu=1}^2\frac{g'}{g-\a_\nu}.$$
Then $\phi$ has only simple poles, exactly at those $\a_3$-points
and poles $(\a_4$-points) that are {\it not} simultaneously simple
for $f$ and $g$, and the usual technique yields
$$T(r,\phi)=\Nquer(r,\a_3)- N_\s(r,\a_3)+\Nquer(r,\a_4)-
 N_\s(r,\a_4)+S(r).$$
Since
 $$\phi(z_\infty)^2=(2\a_3-\a_1-\a_2)^2\Psi(z_\infty)$$
holds at every pole $z_\infty$ that is simple for $f$ and $g$,
$$\begin{array}{rcl}
N_\s(r,\a_4)&\le& 2T(r,\phi)+O(1)\cr &\le& 2\big(\Nquer(r,\a_3)-
 N_\s(r,\a_3)+\Nquer(r,\a_4)- N_\s(r,\a_4)\big)+S(r)\end{array}$$
and $({\rm R}_a)$ for $\kappa=4$ and $\nu=3$ follow, provided
$\phi^2\not\equiv(2\a_3-\a_1-\a_2)^2\Psi$ holds. On the other
hand, $\phi^2=(2\a_3-\a_1-\a_2)^2\Psi$ implies that $\phi$ is a
small function, thus $f$ and $g$ share the values $\a_3$ and
$\a_4=\infty$ by counting multiplicities, and hence all values
$\a_\nu$ by the (2+2)--Theorem. For $2\a_3-\a_1-\a_2=0$
(equivalently $(\a_1,\a_2,\a_3,\infty)=-1$) we obtain the better
inequality
$$2N_\s(r,\a_4)+N_\s(r,\a_3)
 \le\Nquer(r,\a_3)+\Nquer(r,\a_4)+S(r),$$
by counting the zeros of $\phi$ rather than those of
$\phi^2-(2\a_3-\a_1-\a_2)^2\Psi$.

To prove (${\rm R}_b$) we may assume $\kappa=4$, $\a_4=\infty$ and
$\lambda=3$. Then
$$\phi=\frac{f''}{f'}-2\frac{f'}{f-\a_3}+
\frac{f'(f-\a_3)}{(f-\a_1)(f-\a_2)}-
\frac{g''}{g'}-2\frac{g'}{g-\a_3}-\frac{g'(g-\a_3)}{(g-\a_1)(g-\a_2)}$$
is regular at poles of $f$ and $g$ ($\a_4=\infty$), has simple
poles exactly at those $\a_\nu$-points ($1\le\nu\le 3$) of $f$ and
$g$ that have different multiplicities, thus
\be\label{Trphi}T(r,\phi)=N(r,\phi)+S(r)=\sum_{\nu=1}^3\big(\Nquer(r,\a_\nu)-\Nquer_\s(r,\a_\nu)\big)+S(r)\ee
holds. Now $\phi$ vanishes at $\a_3$-points which are
simultaneously simple for $f$ and $g$, hence $\phi\not\equiv 0$
implies
$$N_\s(r,\a_3)\le \sum_{\nu=1}^3\big(\Nquer(r,\a_\nu)-N_\s(r,\a_\nu)\big)+S(r),$$
thus (${\rm R}_b$) for $\kappa=4$ and $\lambda=3$. On the other
hand, the conclusion of the Four-Value-Theorem holds if $\phi$
vanishes identically.

To prove (${\rm R}_c$), for $\kappa=4$ and $\a_4=\infty$, say, we
just consider
 $$\phi=\frac{f'}{(f-\a_1)(f-\a_2)(f-\a_3)}-\frac{g'}{(g-\a_1)(g-\a_2)(g-\a_3)};$$
then $\phi$ has poles exactly at those $\a_\nu$-points
($1\le\nu\le 3$) that have different multiplicities for $f$ and
$g$, and thus satisfies {\it a fortiori} (\ref{Trphi}).
Since $\phi$ vanishes at poles ($\a_4$-points) of $f$ and $g$,
this yields
 $$\Nquer(r,\a_4)\le\sum_{\nu=1}^3\big(\Nquer(r,\a_\nu)-N_\s(r,\a_\nu)\big)+S(r)$$
provided $\phi$ does not vanish identically. If, however, $\phi=
0$, then $f$ and $g$ share the values $\a_\mu$ $(1\le\mu\le 3)$ by
counting multiplicities, and thus the hypothesis of the
(3+1)--Theorem and the conclusion of the Four-Value-Theorem holds.

To prove (R$_d$) we may assume $\a_\mu\in\C$ ($1\le\mu\le 4)$ and
consider
$$\phi=\frac{f'}{(f-\a_1)(f-\a_2)}-\frac{g'}{(g-\a_1)(g-\a_2)}.$$
Then $\phi$ has poles exactly at those $\a_1$- and $\a_2$-points
that are not simultaneously simple for $f$ and $g$, and is regular
at poles of $f$ and of $g$, hence
$$T(r,\phi)=\sum_{\mu=1}^2\big(\Nquer(r,\a_\mu)-N_\s(r,\a_\mu)\big)+S(r)$$
follows. On the other hand,
$\ds\phi(z_\rho)=(-1)^\rho\frac{f'(z_\rho)-g'(z_{\rho})}{\prod\limits_{\mu\ne
 \rho}(\a_\rho-\a_\mu)}\,(\a_4-\a_3)$
holds at any $\a_\rho$-point $z_\rho$ ($\rho=3,4$), and
$\ds\Psi(z_\rho)=\frac{(f'(z_\rho)-g'(z_{\rho}))^2}{\prod\limits_{\mu\ne
 \rho}(\a_\rho-\a_\mu)^2}$
holds whenever $z_\rho$  is simple for $f$ and $g$. Thus if
$\phi^2\not\equiv(\a_4-\a_3)^2\Psi$, the assertion (for
$\kappa=3$, $\nu=4$) follows from the First Main Theorem of
Nevanlinna:
$$\begin{array}{rcl}
\ds\sum_{\mu=3}^4N_\s(r,\a_\mu)&\le&\ds
T\Big(r\frac1{\phi^2-(\a_4-\a_3)^2\Psi}\Big)+O(1)=2T(r,\phi)+S(r)\cr
&\le& \ds 2\sum_{\mu=1}^2\big(\Nquer(r,\a_\mu)-
N_\s(r,\a_\mu)\big)+S(r).\end{array}$$ If, however,
$\phi^2=(\a_4-\a_3)^2\Psi$, then $f$ and $g$ share the values
$\a_1$ and $\a_2$ by counting multiplicities, hence the conclusion
of the Four-Value-Theorem holds.

Finally, (R$_e$) and (R$_f$) follow by adding up inequality
(R$_b$) for $\kappa\ne\lambda$ and (R$_c$)
for $\kappa=1,\ldots,4$, respectively. \Ende

\section{\bf Towards or way off a (1+3)--Theorem?}

 Assuming that $f$ and $g$ share the values $\a_\nu\in\C$
$(1\le \nu\le 3)$ and $\a_4=\infty$, we set
$$\begin{array}{rcl}\ds\phi_{f}=\frac{f'}{(f-\a_1)(f-\a_2)(f-\a_3)}&{\rm
and}&\ds\phi_{g}=\frac{g'}{(g-\a_1)(g-\a_2)(g-\a_3)},\cr
\Phi_{f}=(f-g)\phi_{f}&{\rm and}& \Phi_{g}=(f-g)\phi_{g},\cr
\Phi=\Phi_{f}/\Phi_{g}=\phi_{f}/\phi_{g}&{\rm
and}&\Psi=\Phi_{f}\Phi_{g}
~\textrm{(Mues'~function)}.\end{array}$$

{\it Example.}\quad\begin{tabular}{|r|c|c|c|c|}\hline functions
$f,$ $g$&$\Phi_{f}$&$\Phi_{g}$&$\Psi$&$\Phi$\cr\hline
 $e^z,$ $e^{-z}$&$e^{-z}$&$e^z$&$1$&$e^{-2z}$\cr\hline
 Gundersen's &$1-e^z$&$\ds\frac{8}{1-e^z}$&$8$&$\ds\frac{(1-e^z)^2}8$\cr\hline
 Reinders'&$\ds\frac{12\sqrt{3}}{u+1}$&$\ds\frac{12(u+1)}{\sqrt{3}}$&$144$&$\ds\frac{9}{(u+1)^2}$
 \cr\hline\end{tabular}

\bigskip
{\bf Key Observations.}

\begin{itemize}\item {\it If the value $\a_4=\infty$ is shared by counting
multipli\-cities, then $\Phi_f$ and $\Phi_g$ are entire functions
satisfying $N(r,1/\Phi_f)+N(r,1/\Phi_g)=S(r)$. \item If $T(r)$ has
finite lower order $\ds\liminf_{r\to\infty}\frac{\log T(r)}{\log
r},$ then
$$\qquad\Phi_{f}=p_{f}e^{Q},~
\Phi_{g}=p_{g}e^{-Q},~\Psi=p_fp_g,~\Phi=\frac{p_f}{p_g}\,e^{2Q},
 {\rm~and~}T(r)\asymp r^{\deg Q}$$
hold with polynomials $p_f,$ $p_g,$ and $Q$.}\end{itemize}

\medskip{\it Proof.}  Noting that
$$\Psi=\frac{f^2f'}{P(f)}\frac{g'}{P(g)}+2\frac{ff'}{P(f)}\frac{gg'}{P(g)}+
 \frac{f'}{P(f)}\frac{g^2g'}{P(g)}$$
(with $P(w)=\prod\limits_{\nu=1}^3(w-\a_\nu)$ or
$P(w)=\prod\limits_{\nu=1}^4(w-\a_\nu)$) holds, the lemma on the
logarithmic derivative gives
 $$T(r_k,\Psi)=m(r_k,\Psi)=O(\log r_k)$$
on some sequence $r_k\to\infty$. Hence $\Psi$ is a polynomial,
$\Phi_{f}$ and $\Phi_{g}$ have only finitely many zeros and
satisfy
 $$\log T(r_k,\Phi_f)+\log T(r_k,\Phi_g)=O(\log r_k).$$
Thus the assertion on $\Phi_f$, $\Phi_g$ and $\Psi,$ holds, since
entire functions $e^{h(z)}$ have finite lower order if and only if
$h$ is a polynomial (by the Borel-Carath\'eodory inequality or the
lemma on the logarithmic derivative). The assertion on $T(r)$
follows from the subsequent theorem and $T(r,e^Q)\asymp r^{\deg Q}$.\Ende

\medskip{\bf Theorem.} {\it If $f$ and $g$ share four values, at
least one of them by counting multiplicities, then either the
conclusion of the Four-Value-Theorem or else
 $${\ts\frac5{209}}T(r)\le T(r,\Phi)+S(r)\le 2T(r)-2\Nquer(r,\infty)+S(r)$$
holds.}\medskip

{\it Proof.} We assume that $\a_4=\infty$ is shared by counting
multiplicities, and first suppose that
$\Phi=\Phi_f/\Phi_g=\phi_f/\phi_g$ is non-constant. Then from
$$\begin{array}{rcl}
N(r,1/\Phi_g)&\le&N(r,1/\Psi)+S(r)=S(r),\cr
N(r,1/\phi_g)&=&\Nquer(r,\infty)+S(r)\cr
 m(r,\phi_g)&=&S(r)\end{array}$$
(and the same for $g$ replaced by $f$) follows the upper estimate
$$\begin{array}{rcl}
T(r,\Phi)=m(r,\Phi)+S(r)&\le& m(r,1/\phi_{g})+S(r)\cr &=&\ds
N(r,\phi_{g})-N(r,1/\phi_{g})+S(r)\cr &=&\ds\sum_{\nu=1}^3
\Nquer(r,\a_\nu)-\Nquer(r,\infty)+S(r)\cr
 &=&2T(r)-2\Nquer(r,\infty)+S(r).\end{array}$$
Now suppose that $f(z_0)=g(z_0)=\a_\nu$ holds with multiplicities
$\ell_f$ and $\ell_g$, say, hence $\Phi(z_0)=\ell_f/\ell_g$.
Noting that the sequence of $\a_\nu$-points with
$\min\{\ell_f,\ell_g\}>1$ has counting function $S(r)$, and
restricting $\ell_f$ and $\ell_g$ to the range
$\{1,\ldots,\ell\}$, we obtain by Nevanlinna's second main theorem
$$\begin{array}{rcl}
 T(r)&\le&\ds\sum_{\nu=1}^3\Nquer(r,\a_\nu)+S(r)\cr
 &\le&\ds (2\ell+1)T(r,\Phi)+\frac 1{\ell+2}\sum_{\nu=1}^3
\Big[N\Big(r,\frac{1}{f-\a_\nu}\Big)+N(r,\frac
1{g-\a_\nu}\Big)\Big]+S(r)\cr
 &\le&\ds (2\ell+1)T(r,\Phi)+\frac 6{\ell+2}T(r)+S(r),\end{array}$$
hence
 \be\label{defekt}T(r)\le\frac{(2\ell+1)(\ell+2)}{\ell-4}T(r,\Phi)+S(r)\quad(\ell>4)\ee
holds. The factor $\frac{209}{5}$ is obtained for $\ell=9$. If,
however, $\Phi$ is constant, then we have actually
$\Phi\equiv\ell$ or $\Phi\equiv1/\ell$ for some $\ell\in\N$. This
yields, in the first case, say,
 $$\prod_{\nu=1}^3(f-\a_\nu)=C\prod_{\nu=1}^3(g-\a_\nu)^\ell\quad(C\ne 0 {\rm ~some~constant}).$$
From $T(r,f)\sim \ell T(r,g)$ then follows $\ell=1$, hence $f$ and
$g$ share all values $\a_\nu$ by counting multiplicities. \Ende

Combining the previous result with the second Key Observation, we
obtain:

\medskip{\bf Corollary (Yi \& Li \cite{yili}).} {\it Suppose $f$ and
$g$ share four values $\a_\nu$, at least one of them by counting
multiplicities. Then either both functions have infinite lower
order or else have equal finite integer order that also equals the
lower order.}\medskip

{\it Remark.} The authors of \cite{yili} proved in a way similar
to ours the inequality
 $${\ts\frac1{77}}T(r)\le T(r,\Phi)+S(r)\le 4T(r)+S(r).$$

Each of the counterexamples (by Gundersen, Reinders and the
author) are (can be reduced to) either rational functions of $e^z$
or else elliptic functions. Since elliptic functions have no
deficient value, we obtain:

\medskip{\bf Corollary.} {\it No pair of elliptic functions can share
four values, at least one of them by counting multiplicities.}

\medskip{\bf Theorem.} {\it Let $f$ and $g$ be $2\pi i$-periodic functions of finite
order sharing four values, one of them by counting multiplicities.
Then $f$ and $g$ are rational functions of $e^{z}$.}

\medskip{\it Proof.} The functions $\Psi=\Phi_f\Phi_g,$ $\Phi_f$ and
$\Phi_g$ are entire of finite order and also $2\pi i$-periodic.
But $\Psi$ is a polynomial, hence constant, and thus $\Phi_f$ and
$\Phi_g$ are zero-free. This yields
 \be\label{ehochmz}\Phi_f(z)=e^{mz+c_f}{\rm ~and~} \Phi_g(z)=e^{-mz+c_g}\ee
for some $m\in\Z\setminus\{0\}$ and complex constants $c_f$ and
$c_g$. From $T(r,\Phi_f/\Phi_g)\sim\frac{2|m|}\pi r$, hence
$T(r)\asymp r$ then follows that $f$ and $g$ are rational
functions of $e^z$. \Ende

\medskip{\it Remark.} Generally spoken, meromorphic functions $h(z)=R(e^z)$,
where $R$ is rational with $\deg R=d>1$, have Nevanlinna
characteristic $T(r,h)\sim\frac d\pi r$ and deficient values
$R(0)$ and $R(\infty)$, which may, of course, coincide. If, e.g.,
$R(u)\sim \a+bu^{-\rho}$ as $u\to\infty$, then the contribution of
the right half plane to $m\big(r,1/(f-\a))$ is $\sim\frac\rho\pi
r$. In Gundersen's example we have $\delta(\a,f)=1/2$ for $\a=0,
1,$ and $\delta(\a,g)=1/2$ for $\a=\infty,-1/8.$

Let $f(z)=R(e^z)$ and $g(z)=S(e^z)$ share four values $\a_\nu$
($\in\C$ for technical reasons), without assuming anything about
multiplicities. Then Mues' function $\Psi=\Phi_f\Phi_g$ is a
non-zero constant, and from $R(u)\sim \a_\mu+bu^{-\rho}$ and
$S(u)\sim \a_\nu+cu^{-\sigma}$ $(bc\ne 0,$ $\rho,\sigma>0)$ as
$u\to\infty$ follows
 $$\begin{array}{rcl}\Phi_f(z)&\sim& b_1(\a_\mu-\a_\nu)+O(e^{-\min\{\rho,\sigma\}\Re z})\quad(b_1\ne 0)\cr
 \Phi_f(z)&\sim& c_1(\a_\mu-\a_\nu)+O(e^{-\min\{\rho,\sigma\}\Re z}) \quad(c_1\ne 0)\end{array}
 \quad{\rm as~}\Re z\to +\infty.$$
Thus $\a_\mu\ne \a_\nu$, and a similar result near $u=0$ shows
that at least one of the values $\a_\nu$ (but no value $b\ne
\a_\nu$) is deficient for $f$, and the same is true for $g$ (and
the same or some other value $\a_\mu$). If $\a_4,$ say, is shared
by counting multiplicities, then
$\delta(\a_4,f)=\delta(\a_4,g)>0.$ More precisely we have
$R(u)\sim \a_4+bu^{-\rho}$ $(u\to\infty),$ $S(u)\sim
\a_4+cu^{\rho}$ $(u\to 0),$
$\delta(\a_4,f)=\delta(\a_4,g)=\rho/d$, and $R(0)=\a_\mu$,
$S(\infty)=\a_\nu$ for some $\mu,\nu<4$. We switch now to values
$\a_1, \a_2, \a_3$ and $\a_4=\infty$ and assume that $\a_4$ is
shared by counting multiplicities, and also that
$R(\infty)=S(0)=\infty$ holds. From $m(r,f)\sim
m(r,g)\sim\frac\rho\pi r$ for some $\rho$ $(1\le\rho<d)$ we obtain
 $$R(u)=\frac{P(u)}{Q(u)},~ S(u)=\frac{\tilde P(u)}{u^\rho
 Q(u)},~\deg \tilde P\le d {\rm~and~} \rho+\deg Q=\deg P=d;$$
the zeros of $Q$ are simple and $\ne 0.$

\section{\bf Functions of finite order} 

Under certain circumstances it may happen that a problem for
meromorphic function of arbitrary order of growth can be reduced
to a problem for functions of finite order by the well-known
Zalcman Lemma \cite{zalcman}. A prominent example can be found in
\cite{bewere}. This might also be the case here, although there
are some obstacles, as will be seen later.\medskip

{\bf Rescaling Lemma.} {\it Let $f$ and $g$ share the values
$\a_\nu$ $(1\le\nu\le 4)$. Then either the spherical derivatives
$f^\#$ and $g^\#$ are bounded on $\C$, or else there exist
sequences $(z_k)\subset\C$ and $(\rho_k)\subset(0,\infty)$ with
$\rho_k\to 0$, such that the sequences $(f_k)$ and $(g_k)$,
defined by
$$ f_k(z)=f(z_k+\rho_kz){\rm ~and~} g_k(z)=g(z_k+\rho_kz),$$
simultaneously tend to non-constant meromorphic 
functions $\hat f$ and $\hat g$, respectively; $\hat f$ and $\hat
g$ share the values $\a_\nu$ and have bounded spherical
derivatives.}\medskip

{\it Proof.} We assume that $f^\#$ is not bounded. Then the
existence of the sequences $(z_k)$ (tending to infinity) and
$(\rho_k)$ can be taken for granted for the function $f$ by
Zalcman's Lemma.

If we assume that the sequence $(g_k)$ is not normal on $\C$, then
again by Zalcman's Lemma there exist sequences
$(k_\ell)\subset\N$, $\hat z_\ell\to z_0\in\C$ and
$\sigma_\ell\downarrow 0$, such that $\hat g_\ell=g_{k_\ell}(\hat
z_\ell+\sigma_\ell z)$ tends to some non-constant meromorphic
function $\tilde g$. Then on one hand, the corresponding sequence
$(\hat f_\ell)$ tends to a constant, while on the other hand every
limit function of $(\hat f_\ell)$ shares the values $\a_\nu$ with
$\tilde g$ by Hurwitz' Theorem. Since $\tilde f$ assumes at least
two of these values, our assumption on the sequence $(g_k)$ was
invalid. We may assume that $(g_k)$ tends to $\hat g$. Then $\hat
g$ shares the values $\a_\nu$ with $\hat f$, hence is
non-constant.

Finally, if we assume that $\hat g^\#$ is unbounded on $\C$, then
we may apply the first argument to the functions $\hat g$ and
$\hat f$ in this order, but now with the {\it a priori} knowledge
that $\hat f^\#$ is bounded. Then some sequence $\hat g(\tilde z_k
+\tau_k z)$ tends to some non-constant limit, while a sub-sequence
of the sequence $\hat f(\tilde z_k +\tau_k z)$ tends to a constant
(the sequence of spherical derivatives is $O(\tau_k)$). This
proves the Rescaling Lemma. \Ende

\medskip{\it Remark.} Till now, however, it cannot be excluded that
\begin{itemize}\item[(a)] a shared value gets lost in the sense that
it becomes a Picard value for $\hat f$ and $\hat g$, although it
is not for $f$ and $g$; \item[(b)] multiplicities get lost in the
sense that $\hat f$ and $\hat g$ share some value $\a_\nu$ by
counting multiplicities, although $f$ and $g$ do not; \item[(c)]
$\hat f=\hat g$ (worst case).
\end{itemize}
In the third case everything is lost. But even if $\hat f\ne \hat
g$ holds and the (3+1)-Conjecture turns out to be true for
functions with bounded spherical derivative, we can only deduce
$(\a_1,\a_2,\a_3,\a_4)=-1$, except if we are able to rule out also
the first and second case. Nevertheless we proceed in this
direction.

\medskip{\bf Theorem.} {\it Suppose that $f$ and $g$
share the values $\a_\nu$ $(1\le\nu\le 4)$ and have bounded spherical derivative. Then Mues' function
$\Psi$ is a constant.
Moreover, if one value is shared by counting multiplicities, then
 \be\label{PhifPhig}\Phi_f(z)=e^{Q(z)+c_f} {\rm ~and~}
 \Phi_g(z)=e^{-Q(z)+c_g}\ee
holds, with $Q$ some polynomial of degree one or two.}\medskip

{\it Proof.} From $f^\#+g^\#\le C$ follows that $f$ and $g$ have
finite order (actually the order is at most two), and by
Nevanlinna's lemma on the logarithmic derivative, $\Psi$ is a
polynomial. We assume $\Psi(z)\sim c z^m$ as $z\to\infty$ for some
$m\ge 1$ and $c\ne 0$, and consider, for some sequence
$z_n\to\infty$, the functions
$$ f_n(z)=
 f(z_n+z_n^{-m/2}z){\rm ~and~} g_n(z)= g(z_n+z_n^{-m/2}z),$$
which also share the values $\a_\nu$ and have Mues' function
\be\label{limespsi}\Psi_n(z)=z_n^{-m}\Psi(z_n+z_n^{-m/2}z)\sim
 c\quad(n\to\infty),\ee
while obviously $f_n$ and $g_n$ tend to constants. If the sequence
$(z_n)$, which is quite arbitrary, can be chosen in such a way
that
 $$f(z_n)\to b\ne\a_\nu {\rm ~and~}g(z_n)\to b'\ne\a_\nu\quad(1\le \nu\le
 4),$$
then $f_n$ and $g_n$ tend to constants $b$ and $b'$, hence
$\Psi_n$ tends to $0$ in contrast to relation (\ref{limespsi}).
Now for any sequence $(\hat z_n)$ the sequences $\hat
f_n(z)=f(\hat z_n+z)$ and $\hat g_n(z)=g(\hat z_n+z)$ are normal
on $\C$, we may assume that (some sub-sequence of $\hat f_n$,
again denoted by) $\hat f_n$ tends to some {\it non-constant}
limit function $\hat f$, e.g.\ by choosing $\hat z_n\to\infty$
such that $|\hat f'(\hat z_n)|=1$. Since $\hat g_{n}$ shares the
values $\a_\nu$ with $\hat f_{n}$ we may (by normality and
Picard's theorem) also assume that $\hat g_n\to \hat g\not\equiv
const.$ Now we choose $z_0$ such that $\hat f(z_0)=b\ne\a_\nu$ and
$\hat g(z_0)=b'\ne\a_\nu$, and set $z_n=\hat z_n+z_0$.


Finally, since $\Psi=\Phi_f\Phi_g$ is constant, the entire
functions $\Phi_f$ and $\Phi_g$ are zero-free and have order one or two
(as do $f$ and $g$), hence (\ref{PhifPhig}) holds with
$\deg Q=1$ or $2$. \Ende

\end{document}